\documentclass{article}
\usepackage[utf8]{inputenc}
\usepackage{amsmath}
\usepackage[T1]{fontenc}    
\usepackage{url}            
\usepackage{booktabs}       
\usepackage{amsfonts}       
\usepackage{nicefrac}       
\usepackage{microtype}      
\usepackage{lipsum}
\usepackage{graphicx}
\usepackage[english]{babel}
\usepackage{amsthm}
\usepackage{tabularx}
\usepackage{array}
\usepackage{mathrsfs}
\usepackage{mathtools}
\usepackage{amssymb}

\newtheorem{theorem}{Theorem}[section]
\newtheorem{conjecture}{Conjecture}[section]
\newtheorem{remark}{Remark}[section]
\newtheorem{corollary}{Corollary}[theorem]
\newtheorem{lemma}[theorem]{Lemma}

\usepackage{listings}
\usepackage{xcolor}
\title{A new approach to odd perfect numbers via GCDs}
\author{Jose Arnaldo Bebita Dris \\ 
  \texttt{josearnaldobdris@gmail.com} \\}

\begin{document}

\maketitle

\begin{abstract}
Let $q^k n^2$ be an odd perfect number with special prime $q$.  Define the GCDs
$$G = \gcd\bigg(\sigma(q^k),\sigma(n^2)\bigg)$$
$$H = \gcd\bigg(n^2,\sigma(n^2)\bigg)$$
and
$$I = \gcd\bigg(n,\sigma(n^2)\bigg).$$
We prove that $G \times H = I^2$.  (Note that it is trivial to show that $G \mid I$ and $I \mid H$ both hold.)  We then compute expressions for $G, H,$ and $I$ in terms of $\sigma(q^k)/2, n,$ and $\gcd\bigg(\sigma(q^k)/2,n\bigg)$.  Afterwards, we prove that if $G = H = I$, then $\sigma(q^k)/2$ is not squarefree.  Other natural and related results are derived further.  Lastly, we conjecture that the set
$$\mathscr{A} = \{m : \gcd(m,\sigma(m^2))=\gcd(m^2,\sigma(m^2))\}$$
has asymptotic density zero.
\end{abstract}

\section{Introduction}\label{Sec1}
In what follows, we let $\sigma(x)$ denote the sum of divisors of the positive integer $x$.  We will denote the deficiency of $x$ by $D(x)=2x-\sigma(x)$, the aliquot sum of $x$ by $s(x)=\sigma(x)-x$, and the abundancy index of $x$ by $I(x)=\sigma(x)/x$.

A number $M$ satisfying $\sigma(M)=2M$ is called a \emph{perfect number}.  For example, $6$ and $28$ are perfect since
$$\sigma(6) = 1 + 2 + 3 + 6 = 2 \cdot 6$$
$$\sigma(28) = 1 + 2 + 4 + 7 + 14 + 28 = 2 \cdot 28.$$
The Euclid-Euler Theorem states that $M$ is an even perfect number if and only if
$$M = (2^t - 1){2^{t-1}},$$
where $2^t - 1$ (and therefore $t$) is prime.  (If $t$ is prime, then $2^t - 1$ is not necessarily prime.)  Primes of the form $2^t - 1$ are called Mersenne primes.  Currently, there are $51$ known Mersenne primes (with the latest being discovered by the Great Internet Mersenne Prime Search in December of 2018), corresponding to $51$ even perfect numbers \cite{GIMPS}.

It is currently unknown whether there are infinitely many even perfect numbers.  It has been conjectured, and is widely believed, that no odd perfect numbers exist.  (There are no odd perfect numbers less than ${10}^{1500}$ \cite{OchemRao}, making the existence of an odd perfect number appear very unlikely.)

Euler proved that an odd perfect number $N$ must necessarily have the so-called \emph{Eulerian form}
$$N = q^k n^2$$
where $q$ is the special prime satisfying $q \equiv k \equiv 1 \pmod 4$ and $\gcd(q,n)=1$.  Descartes, Frenicle, and subsequently Sorli conjectured that $k=\nu_{q}(N)=1$ always holds \cite{Beasley}.  Sorli predicted that $k=1$ after testing large numbers with eight distinct prime factors for perfection \cite{Sorli}.

Dris conjectured that $q^k < n$ \cite{Dris2}, on the basis of the result $I(q^k) < \sqrt[3]{2} < I(n)$.

We state these conjectures here for ease of reference later on.

\begin{conjecture}\label{DFS}
If $q^k n^2$ is an odd perfect number given in Eulerian form, then $k=1$.
\end{conjecture}

\begin{conjecture}\label{Dris}
If $q^k n^2$ is an odd perfect number given in Eulerian form, then $q^k < n$.
\end{conjecture}

Dris \cite{Dris} showed that the equation
$$i(q)=\frac{\sigma(n^2)}{q^k}=\frac{2n^2}{\sigma(q^k)}=\frac{D(n^2)}{s(q^k)}=\frac{2s(n^2)}{D(q^k)}=\gcd(n^2,\sigma(n^2))$$
holds. We also know that the index $i(q)$ is an integer which is at least $3$ by a result of Dris \cite{Dris2}.  (The lower bound on $i(q)$ has since been improved by several authors.)

Furthermore, we can express $i(q)$ as
$$i(q)=q\sigma(n^2) - 2(q-1)n^2.$$

Set $E=n$, $F=\sigma(q^k)/2$, and $K = \gcd(E,F)$.

In this note, we compute expressions for the following GCDs
$$G = \gcd\bigg(\sigma(q^k),\sigma(n^2)\bigg)$$
$$H = \gcd\bigg(n^2,\sigma(n^2)\bigg)$$
and
$$I = \gcd\bigg(n,\sigma(n^2)\bigg).$$
It turns out that it is possible to express all of them in terms of $E$, $F$, and $\gcd(E,F)$.

As far as the author is aware, the approach presented in this paper is new and has not been considered before in the literature.

\section{Preliminaries}\label{Sec2}
Define
$$G = \gcd\bigg(\sigma(q^k),\sigma(n^2)\bigg)$$
$$H = \gcd\bigg(n^2,\sigma(n^2)\bigg)$$
and
$$I = \gcd\bigg(n,\sigma(n^2)\bigg).$$

The following lemma gives an identity that relates the values of $G, H,$ and $I$.
 
\begin{lemma}\label{G times H equals I squared}
If $N = q^k n^2$ is an odd perfect number given in Eulerian form, then $G \times H = I^2$.
\end{lemma}

\begin{proof}
We have
$$\sigma(q^k) = \frac{2n^2}{i(q)}$$
and
$$\sigma(n^2) = {q^k}{i(q)},$$
so that we get
$$G=\gcd\left(\sigma(q^k),\sigma(n^2)\right) = \gcd\bigg(\frac{2n^2}{i(q)}, {q^k}{i(q)}\bigg) = \frac{\gcd\Bigg(n^2, \bigg(i(q)\bigg)^2\Bigg)}{i(q)} = \frac{\Bigg(\gcd\bigg(n, i(q)\bigg)\Bigg)^2}{i(q)}$$
$$= \frac{\Bigg(\gcd\bigg(n, \gcd(n^2, \sigma(n^2))\bigg)\Bigg)^2}{i(q)} = \frac{\Bigg(\gcd\bigg(\gcd(n, n^2),\sigma(n^2)\bigg)\Bigg)^2}{i(q)} = \frac{\Bigg(\gcd\bigg(n,\sigma(n^2)\bigg)\Bigg)^2}{i(q)} = \frac{I^2}{H},$$
since $\sigma(n^2)$ is odd and $\gcd(q,n)=1$.
\end{proof}

Using the identity in Lemma \ref{G times H equals I squared}, we can now derive the following divisibility conditions.

\begin{lemma}\label{G divides I and I divides H}
Suppose that $N = q^k n^2$ is an odd perfect number given in Eulerian form.  Then $G$ divides $I$ and $I$ divides $H$.
\end{lemma}

\begin{proof}
The proof of the divisibility constraint $I \mid H$ follows from the GCD property
$$\bigg(a \mid b\bigg) \implies \bigg(\gcd(a,c) \mid \gcd(b,c)\bigg).$$
Afterwards, the proof of the divisibility constraint $G \mid I$ then follows from Lemma \ref{G times H equals I squared}.
\end{proof}

Set 
$$J=\frac{I}{G}=\frac{H}{I}.$$

By Lemma \ref{G times H equals I squared} and Lemma \ref{G divides I and I divides H}, $J$ is an (odd) integer.

The following lemma computes the value of $J$, in terms of $E, F,$ and $\gcd(E,F)$.  (The proof is due to the anonymous MSE user mathlove \cite{mathlove}.)

\begin{lemma}\label{Value for J}
If $N = q^k n^2$ is an odd perfect number given in Eulerian form, then we obtain
$$J = \frac{n}{\gcd\bigg(\sigma(q^k)/2,n\bigg)}.$$
\end{lemma}

\begin{proof}
We have
$$H = \frac{n^2}{\sigma(q^k)/2}.$$
Hence, we obtain
$$J = \frac{H}{I} = \frac{n^2}{\sigma(q^k)/2\cdot\gcd\bigg(n,\sigma(n^2)\bigg)}=\frac{n^2}{\gcd\bigg(n\cdot\sigma(q^k)/2,\sigma(q^k)\sigma(n^2)/2\bigg)}=\frac{n^2}{\gcd\bigg(n\cdot\sigma(q^k)/2,q^k n^2\bigg)}$$
$$=\dfrac{n}{\gcd\bigg(\sigma(q^k)/2,q^k n\bigg)}=\frac{n}{\gcd\bigg(\sigma(q^k)/2,n\bigg)}.$$
\end{proof}

\begin{remark}\label{General Equation for G, H, I, and J}
Notice that we then have
$$J = \frac{I}{G} = \frac{H}{I} = \sqrt{\frac{H}{G}}$$
so that
$$H = G \times {J^2}.$$
\end{remark}

We can now compute expressions for $I$ and $G$, using Lemma \ref{Value for J}.

\begin{lemma}\label{Values for I and G}
If $N = q^k n^2$ is an odd perfect number given in Eulerian form, then we obtain
$$I = \Bigg(\frac{n}{\sigma(q^k)/2}\Bigg)\cdot{\gcd\bigg(\sigma(q^k)/2,n\bigg)}$$
and
$$G = \frac{\Bigg(\gcd\bigg(\sigma(q^k)/2,n\bigg)\Bigg)^2}{\sigma(q^k)/2}.$$
\end{lemma}

\begin{proof}
The proof is trivial.
\end{proof}

\section{What happens when $J=1$?}\label{Sec4}
Let us examine the case $J=1$ to see whether it is interesting.

First, we prove the following unconditional lemma.

\begin{lemma}\label{When J is 1}
Suppose that $N = q^k n^2$ is an odd perfect number given in Eulerian form.  Then $J = 1$ holds if and only if $n \mid \sigma(q^k)/2$.
\end{lemma}

\begin{proof}
Recall that from Lemma \ref{Value for J}, we have 
$$J = \frac{n}{\gcd\bigg(\sigma(q^k)/2,n\bigg)}.$$
This is equal to one if and only if
$$\gcd\bigg(\sigma(q^k)/2,n\bigg) = n,$$
which holds if and only if $n \mid \sigma(q^k)/2$.
\end{proof}

\begin{remark}\label{Remark on J equal to 1}
Note that $J = 1$ if and only if $G = H = I$ holds.
\end{remark}

If $J=1$, then from Lemma \ref{When J is 1}, we get $n \mid \sigma(q^k)/2$, from which we obtain $n < q^k$.  But Brown \cite{Brown} proved the estimate $q < n$ in 2016.  Hence, $J = 1$ implies that $k > 1$.  We record this in the succeeding proposition.

\begin{lemma}\label{When J is 1 then k is not 1}
Suppose that $N = q^k n^2$ is an odd perfect number given in Eulerian form.  If $J = 1$, then both Conjecture \ref{DFS} and Conjecture \ref{Dris} are false.
\end{lemma}

Recall from Remark \ref{General Equation for G, H, I, and J} that $H = G \times J^2$.  Since $H \geq 3$ holds \cite{Dris2}, $G = J = 1$ is not true.  We therefore get the following proposition.

\begin{lemma}\label{When J is 1 then G is not 1}
Suppose that $N = q^k n^2$ is an odd perfect number given in Eulerian form.  If $J = 1$, then $G \neq 1$ holds.
\end{lemma}

\begin{theorem}\label{Equivalent Conditions 2}
Suppose that $N = q^k n^2$ is an odd perfect number given in Eulerian form.  The following conditions are equivalent to $J=1$:
\begin{enumerate}
{
\item{$n \mid \sigma(q^k)/2$}
\item{$\sigma(n^2) \mid q^k n$.}
}
\end{enumerate}
\end{theorem}

\begin{proof}
The proof follows from Lemma \ref{When J is 1}, and by writing the equation
$$\frac{\sigma(n^2)}{n} = \frac{q^k n}{\sigma(q^k)/2}$$
in the form
$$\frac{q^k n}{\sigma(n^2)} = \frac{\sigma(q^k)/2}{n}.$$
\end{proof}

\section{What happens when $F$ is squarefree?}\label{Sec3}
We rewrite the equation
$$\frac{\sigma(n^2)}{q^k} = \frac{n^2}{\sigma(q^k)/2}$$
in the form
$$\frac{\sigma(n^2)}{n} = \frac{q^k n}{\sigma(q^k)/2}$$
to get the succeeding proposition.

The following theorem is similar in spirit to Theorem \ref{Equivalent Conditions 2}.

\begin{theorem}\label{Equivalent Conditions 1}
If $N = q^k n^2$ is an odd perfect number given in Eulerian form, then the following conditions are equivalent:
\begin{enumerate}
{
\item{$\sigma(q^k)/2 \mid n$}
\item{$n \mid \sigma(n^2)$}
\item{$G = \sigma(q^k)/2$}
\item{$I = n$}
}
\end{enumerate}
\end{theorem}

\begin{proof}
The equivalence of the first two conditions follows from the fact that $\gcd(q^k,\sigma(q^k))=1$.

Next, we show that the third condition is equivalent to the first.  Recall from Lemma \ref{Values for I and G} that 
$$G = \frac{\Bigg(\gcd\bigg(\sigma(q^k)/2,n\bigg)\Bigg)^2}{\sigma(q^k)/2}.$$
We then see that $G = \sigma(q^k)/2$ if and only if $\sigma(q^k)/2 \mid n$.

Lastly, we show that the fourth condition is equivalent to the second.  To this end, suppose that
$$n = I = \gcd(n, \sigma(n^2)).$$
By the definition of GCD, it follows that $n \mid \sigma(n^2)$.  Conversely, assume that $n \mid \sigma(n^2)$.  Then we obtain
$$I = \gcd(n, \sigma(n^2)) = n$$
by the definition of GCD, and we are done.
\end{proof}

Under the condition that $\sigma(q^k)/2$ is squarefree, we get the following result.

\begin{theorem}\label{When F is squarefree}
Suppose that $N = q^k n^2$ is an odd perfect number given in Eulerian form.  If $F = \sigma(q^k)/2$ is squarefree, then $J \neq 1$.
\end{theorem}

\begin{proof}
Suppose that $\sigma(q^k)/2$ is squarefree.  (Note that, since $\sigma(q^k)/2 \mid n^2$ holds in general, then this hypothesis implies that $\sigma(q^k)/2 \mid n$ is true.) Assume to the contrary that $J = 1$.  By Theorem \ref{Equivalent Conditions 2}, we have $n \mid \sigma(q^k)/2$.  Since $\sigma(q^k)/2$ and $n$ are both positive, we obtain $\sigma(q^k)/2 = n$.  This contradicts a result of Steuerwald in 1937 \cite{Steuerwald}, who proved that $n$ must contain a square factor.
\end{proof}

\subsection{What happens when $H$ is squarefree?}\label{Subsec3.1}
Suppose that $H$ is squarefree.  By Lemma \ref{G times H equals I squared}, we have the equation
$$G \times H = I^2.$$
Since this implies $H \mid I^2$, it follows that $H \mid I$. However, by Lemma \ref{G divides I and I divides H}, we have $I \mid H$.  Since $I$ and $H$ are both positive, we obtain $I = H$.  This means that $J = 1$.  By the contrapositive to Theorem \ref{When F is squarefree}, we finally have $\sigma(q^k)/2$ is not squarefree.

We record the immediately preceding results in the following propositions.

\begin{theorem}\label{When H is squarefree thm}
Suppose that $N = q^k n^2$ is an odd perfect number given in Eulerian form.  If $H = G \times J^2$ is squarefree, then $J = 1$.
\end{theorem}

\begin{corollary}\label{When H is squarefree cor}
Suppose that $N = q^k n^2$ is an odd perfect number given in Eulerian form.  If $H$ is squarefree, then $F = \sigma(q^k)/2$ is not squarefree.
\end{corollary}

\section{On the equation $H = I$}\label{Sec5}
Recall that we have the biconditional
$$J = 1 \iff H = I$$
from Remark \ref{Remark on J equal to 1}.

In this section, we shall attempt a naive determination of the asymptotic density of positive integers $m$ satisfying the equation $\gcd(m,\sigma(m^2))=\gcd(m^2,\sigma(m^2))$.  

The author tried searching for examples and counterexamples via Sage Cell Server.

All positive integers from $1$ to $100$ (except for the integer $99$) satisfy the equation.

The following integers in the range $1 \leq m \leq 1000$ \emph{do not} satisfy $\gcd(m,\sigma(m^2))=\gcd(m^2,\sigma(m^2))$.
$$99 = {3^2}\cdot{11}$$
$$154 = 2\cdot 7\cdot 11$$
$$198 = 2\cdot{3^2}\cdot{11}$$
$$273 = 3\cdot 7\cdot 13$$
$$322 = 2\cdot 7\cdot 23$$
$$396 = {2^2}\cdot{3^2}\cdot{11}$$
$$399 = 3\cdot 7\cdot 19$$
$$462 = 2\cdot 3\cdot 7\cdot 11$$
$$469 = 7\cdot 67$$
$$495 = {3^2}\cdot 5\cdot 11$$
$$518 = 2\cdot 7\cdot 37$$
$$546 = 2\cdot 3\cdot 7\cdot 13$$
$$553 = 7\cdot 79$$
$$620 = {2^2}\cdot 5\cdot 31$$
$$651 = 3\cdot 7\cdot 31$$
$$693 = {3^2}\cdot 7\cdot 11$$
$$741 = 3\cdot 13\cdot 19$$
$$742 = 2\cdot 7\cdot 53$$
$$770 = 2\cdot 5\cdot 7\cdot 11$$
$$777 = 3\cdot 7\cdot 37$$
$$792 = {2^3}\cdot{3^2}\cdot 11$$
$$798 = 2\cdot 3\cdot 7\cdot 19$$
$$903 = 3\cdot 7\cdot 43$$
$$938 = 2\cdot 7\cdot 67$$
$$966 = 2\cdot 3\cdot 7\cdot 23$$
$$990 = 2\cdot{3^2}\cdot 5\cdot 11$$

A simple inspection yields that primes and prime powers satisfy the equation, so that there are infinitely many solutions.

The following Pari/GP-routines efficiently determine the numbers and percentages of solutions, up to a certain search limit.  One can easily adjust the range.

\begin{lstlisting}
c=0;for(m=1,10,if(gcd(m,sigma(m^2))==gcd(m^2,sigma(m^2)),c=c+1));print(c,"  ",((c/10)*1.0))

c=0;for(m=1,100,if(gcd(m,sigma(m^2))==gcd(m^2,sigma(m^2)),c=c+1));print(c,"  ",((c/100)*1.0))

\end{lstlisting}

To summarize, we have the table below which shows the counts and percentages of the number of solutions to the equation
$$\gcd(m,\sigma(m^2))=\gcd(m^2,\sigma(m^2)),$$
up to $10$, ${10}^2$, ${10}^3$, ${10}^4$, ${10}^5$, and ${10}^6$, respectively:

\begin{center}
    \begin{tabular}{|c|c|c|}
    \hline
       Upper limit & Count & Percentage \\ \hline
       10& 10 & 100\% \\ \hline
       100& 99 & 99\% \\ \hline
       1000& 974 & 97.4\% \\ \hline
       10000& 9561 & 95.61\% \\ \hline
       100000& 93845 & 93.845\% \\ \hline
       1000000& 923464 & 92.3464\% \\ \hline
    \end{tabular}
\end{center}

The author was only able to test until ${10}^6$ because the Pari/GP interpreter of Sage Cell Server begins to crash as soon as a search limit of ${10}^7$ is specified.

The author thinks this is not a rigorous proof, but it is definitely evidence to suggest that the asymptotic density in question is less than one.

We state and prove this assertion in the following theorem, which the author first conjectured in the year 2020:

\begin{theorem}\label{AsymptoticDensity1}
The asymptotic density $\mathscr{A}$ of positive integers $m$ with
$$\gcd(m,\sigma(m^2))=\gcd(m^2,\sigma(m^2))$$
satisfies
$$\mathscr{A} < 1.$$
\end{theorem}

\begin{proof}
Generalizing the first (counter)example of $99$ is trivial.

If ${3^2}\cdot{11} \parallel m$, then $11 \parallel \gcd(m,\sigma(m^2))$  and  $11^2 \parallel \gcd(m^2,\sigma(m^2))$.  So the asymptotic density in question is less than 
$$1-\frac{2}{3^3}\cdot\frac{10}{11^2} = \frac{3247}{3267} \approx 0.993878.$$

Also, if $3 \parallel m$, then with probability $1$ there exist two distinct primes $y$ and $z$ congruent to $1$ modulo $3$ such that $y \parallel m$ and $z \parallel m$.  In this case, we get $3 \parallel \gcd(m,\sigma(m^2))$  and  $3^2 \parallel \gcd(m^2,\sigma(m^2))$.  So the asymptotic density in question is less than
$$1-\frac{2}{3^2} = \frac{7}{9} \approx 0.\overline{777}.$$
\end{proof}

The real open problem is whether the asymptotic density $\mathscr{A}$ is $0$.  We state this in the succeeding conjecture:

\begin{conjecture}\label{AsymptoticDensity2}
The asymptotic density $\mathscr{A}$ of positive integers $m$ with
$$\gcd(m,\sigma(m^2))=\gcd(m^2,\sigma(m^2))$$
satisfies
$$\mathscr{A} = 0.$$
\end{conjecture}

\begin{remark}\label{Remark Asymptotic Density H equals I}
In an answer to one of the author's questions in MathOverflow, Aaron Meyerowitz \cite{Meyerowitz} (\url{https://mathoverflow.net/users/8008}) made the following assertions regarding Conjecture \ref{AsymptoticDensity2}: "I think the density does go to zero, but quite slowly. If $p \equiv 1 \pmod 6$ is prime then there are two solutions $0<r<s<p-1$ of $$u^2+u+1 \equiv 0 \pmod p.$$ 
If $p\parallel m$ then, with probability $1,$ there are two distinct primes $u$ and $v,$  each congruent to $r \pmod p,$ with $u \parallel m$ and $v \parallel m.$ (Either or both could be congruent to $s$ as well.)  Then $p \parallel \gcd(m,\sigma(m^2))$ while $p^2 \parallel \gcd(m^2,\sigma(m^2)).$ So the asymptotic density for this not to happen is $1-\frac{p-1}{p^2}<1-\frac{1}{p+2}$.  If we can argue that the chance that none of these events happen is asymptotically $\prod(1-\frac{p-1}{p^2})$ over the primes congruent to $1 \pmod 6,$ then that asymptotic density is $0$."
\end{remark}

\section{Some Further Considerations}\label{Sec6}

\subsection{Bounds for $K, G, I$ and $J$}\label{Subsec6.1}
Recall from Section \ref{Sec1} that we have set $E = n$ and $F = \sigma(q^k)/2$.  In this section, we compute a lower bound for $K = \gcd(E,F)$.

We begin with the following proposition.

\begin{theorem}\label{GCD of E and F}
If $N = q^k n^2$ is an odd perfect number given in Eulerian form, then $K \neq 1$.
\end{theorem}

\begin{proof}
Assume to the contrary that $K = 1$.  Then, we can simplify the expression for $J$ (from Lemma \ref{Value for J}) as
$$J = \frac{n}{\gcd\bigg(\sigma(q^k)/2,n\bigg)} = \frac{E}{K} = E.$$
Recall from Remark \ref{General Equation for G, H, I, and J} that
$$J^2 = \frac{H}{G}.$$
We solve for $G$ and then obtain
$$G = \frac{H}{J^2} = \frac{E^2}{F}\cdot\frac{1}{E^2} = \frac{1}{F},$$
whereupon we get a contradiction from $F = \sigma(q^k)/2 \geq 3$ and $G \geq 1$.
\end{proof}

\begin{remark}\label{Remark GCD of E and F}
By Theorem \ref{GCD of E and F}, $K \geq 2$ must hold.  Since $E$ and $F$ are both odd, then we also obtain $K \neq 2$.  Consequently, we have the lower bound
$$K = \gcd(E,F) \geq 3.$$
\end{remark}

As a corollary, we obtain the following (unconditional) bounds for $G$, $I$ and $J$.  (We have also used the definitional property of $K = \gcd(E,F)$, i.e. that $K$ divides both $E$ and $F$.)

\begin{corollary}
If $N = q^k n^2$ is an odd perfect number given in Eulerian form, then the following bounds hold:
\begin{enumerate}
{
\item{$\frac{9}{F} \leq G \leq F$}
\item{$\frac{3E}{F} \leq I \leq E$}
\item{$\frac{E}{F} \leq J \leq \frac{E}{3}$}
}
\end{enumerate}
\end{corollary}

\subsection{On the constraint $\sigma(w^2) \equiv 0 \pmod w$}\label{Subsec6.2}
Lastly, the author also tried checking for examples of numbers $2 \leq w \leq {10}^6$ satisfying the divisibility constraint
$$w \mid \sigma(w^2)$$
using the following Pari-GP script, via Sage Cell Server:
\begin{lstlisting}
for(w=2, 1000000, if((Mod(sigma(w^2),w) == 0),print(w,factor(w))))
\end{lstlisting}

Here is the output:
$$39 = 3 \cdot {13}$$
$$793 = {13} \cdot {61}$$
$$2379 = 3 \cdot {13} \cdot {61}$$
$$7137 = {3^2} \cdot {13} \cdot {61}$$
$$13167 = {3^2} \cdot 7 \cdot {11} \cdot {19}$$
$$76921 = {13} \cdot {61} \cdot {97}$$
$$78507 = {3^2} \cdot {11} \cdot {13} \cdot {61}$$
$$230763 = 3 \cdot {13} \cdot {61} \cdot {97}$$
$$238887 = {3^2} \cdot {11} \cdot {19} \cdot {127}$$
$$549549 = {3^2} \cdot 7 \cdot{11} \cdot {13} \cdot{61}$$
$$692289 = {3^2} \cdot {13} \cdot {61} \cdot {97}$$
$$863577 = {3^2} \cdot {{11}^2} \cdot {13} \cdot {61}$$

Note that all of the known examples are \emph{odd}.  The author double-checked the list of the first $199$ terms of OEIS sequence \textit{A232354} (\url{https://oeis.org/A232354/b232354.txt}) and verified that all of them are odd.  Additionally, all of the terms $w$ in that list \emph{do not satisfy} $\sigma(w^2)/w = d^e$ (where $d$ is prime), except for $w = 39$.

\section{Future Research}\label{Sec7}
We leave the following problem for other researchers to solve.

\begin{conjecture}
If $N = q^k n^2$ is an odd perfect number given in Eulerian form, then unconditionally we have
$$\gcd(\sigma(q^k), \sigma(n^2)) \neq \gcd(n^2, \sigma(n^2)).$$
\end{conjecture}
 
\section{Acknowledgments}\label{Sec8}
The author thanks the anonymous MSE user Peter (\url{https://math.stackexchange.com/u/82961}) for sharing Pari/GP-routines.  The author also would like to give credit to the anonymous MSE user mathlove (\url{https://math.stackexchange.com/users/78967}) for providing a proof of Lemma \ref{Value for J} \cite{mathlove}.  The author is likewise indebted to Aaron Meyerowitz of MathOverflow \cite{Meyerowitz}.

\end{document}